\documentclass[francais]{article-arima}
\author{Horiya Kouadri Habbaze, Ahmed Lakmeche and Abdelkader Lakmeche}
\title{Positive solutions for a second order}
\subtitle{multi-point boundary value problem with delay}
\address{
Biomathematics Laboratory\\
Univ Sidi-Bel-Abb\`es\\
Sidi-Bel-Abb\`es\\
Algeria\\
lakahmed2000@yahoo.fr and lakmeche@yahoo.fr}
\resume{Dans ce travail, nous étudions l'existence des solutions positives pour un problème aux limites en plusieurs points pour une équation différentielle du second ordre avec retard. Sous certaines conditions de croissance sur la non-linéarité, et moyennant le théorème du point fixe de Leray-Schauder, on obtient des conditions suffisantes pour l'existence d'une solution non triviale, ce qui améliorent les résultats de J. Chen {\it et al.} [\ref{y3}].}
\abstract{In this work, we investigate the existence of positive solutions for a multi-point boundary value problem for a second order delay differential equation. Under certain growth conditions on the nonlinearity, and by the mean of Leray-Schauder fixed point theorem, sufficient conditions for the existence of nontrivial solution are obtained, which improve the results of J. Chen {\it et al.} [\ref{y3}].}
\motscles{Solution positive, Equation différentielle à retard, Problème aux limites en plusieurs points, Théorème de point fixe de Leray-Schauder.}
\journal[Rubrique]{\textbf{Revue}}{26}{1}{2016}{1}{11}

\newtheorem{theorem}{Theorem}[section]
\begin{document}
\maketitlepage
\pagestyle{empty}

\section{Introduction}
The boundary value problems for delay differential equations
arise in a variety of areas of applied mathematics, physics and variational
problems of control theory [\ref{y3333}]. Recently, many researchers have
done a great deal of research works upon boundary value problems of lower order differential
equations with delay, and some interesting results were produced, see for example
[\ref{y1}], [\ref{y2}] and [\ref{y5}]-[\ref{y9}].\\
In this work, we study the existence of positive solutions of the following nonlinear multi-point boundary value problem with delay
\begin{eqnarray}
  \begin{array}{l}\label{wq}
   u^{''}(t)+ \lambda a(t)f(t,u(t-\tau))=0,\quad \quad t\in  [0,1],    \\
  u(t)=\beta u(\eta),   \quad \quad  -\tau \leq t \leq 0, \\
u(1)=\alpha u(\eta)    \\
  \end{array}
\end{eqnarray}
 where $0<\tau<1$, $0<\eta<1$,  $0<\alpha<\frac{1}{\eta}$ and $0<\beta<\frac{1-\alpha\eta}{1-\eta}$
 are constants, and $\lambda$ is a positive real parameter.\\
 The paper is organized as follows, in section tow we give definitions and preliminaries, and in section there we give our main results.
\section{Preliminaries} In this section we give some preliminary results.
\begin{definition}\ \\
$u(t)$ is called a positive solution of $(\ref{wq})$ if $u\in \mathrm{C}[-\tau,1]\cap \mathrm{C}^{2}(0,1),$ $u(t)\geq 0$ for $t\in(0,1)$ and satisfies $(\ref{wq}).$
\end{definition}
\begin{lemma}\label{lemm1}\ \\
Let $\beta\neq\frac{1-\alpha\eta}{1-\eta}$. Then for $y\in
C([0, T],\ \mathbf{R})$, the boundary value problem
\begin{eqnarray}
u''(t)+y(t)=0, \quad t\in[0, T],\label{e2.1} \\
u(0)=\beta u(\eta),\quad u(1)=\alpha u(\eta) \label{e2.2}
\end{eqnarray}
has a unique solution
\begin{equation}
 u(t)=\int_0^1 G(t,s)y(s)ds
 \label{e2.3}
\end{equation}
where
\begin{equation}
 G(t,s)= g(t,s)+\frac{\beta+(\alpha-\beta)t }
{(1-\alpha\eta)-\beta(1-\eta)} g(\eta,s)
 \label{e2.3}
\end{equation}
and $$g(t,s)=\left\{
       \begin{array}{ll}
        s(1-t), & 0\leq s\leq t\leq 1, \\
          t(1-s), & 0\leq t\leq s\leq 1.
       \end{array}
     \right.$$
\end{lemma}
\proof{From equation (\ref{e2.1}), we have
$$ u(t)=u(0)+u'(0)t-\int_0^t (t-s)y(s)ds\ :=A+Bt-
\int_0^t (t-s)y(s)ds
$$
with
$$\begin{array}{l}
u(0)=A,\\
u(\eta)=A+B\eta-\int_0^\eta (\eta-s)y(s)ds
\end{array}$$
and
$$u(1)=A+B-\int_0^1 (1-s)y(s)ds.$$
From  $u(0)=\beta u(\eta)$,  we have
$$(1-\beta)A-B\beta\eta=-\beta\int_0^\eta (\eta-s)y(s)ds.$$
From $u(1)=\alpha u(\eta)$, we have
$$(1-\alpha)A+B(1-\alpha\eta)=\int_0^1 (1-s)y(s)ds-\alpha \int_0^\eta (\eta-s)y(s)ds.$$
Therefore,
$$
A=\frac{\beta\eta}{(1-\alpha\eta)-\beta(1-\eta)} \int_0^1
(1-s)y(s)ds  -\frac{\beta }{(1-\alpha\eta)-\beta(1-\eta)}
\int_0^\eta (\eta-s)y(s)ds
$$
and
$$
B=\frac{1-\beta}{(1-\alpha\eta)-\beta(1-\eta)} \int_0^1
(1-s)y(s)ds -\frac{\alpha-\beta} {(1-\alpha\eta)-
\beta(1-\eta)} \int_0^\eta (\eta-s)y(s)ds.
$$
From which it follows that
\begin{eqnarray}
u(t)&=&\frac{\beta\eta}{(1-\alpha\eta)-\beta(1-\eta)} \int_0^1 (1-s)y(s)ds
-\frac{\beta }{(1-\alpha\eta)-\beta(1-\eta)} \int_0^\eta (\eta-s)y(s)ds\nonumber\\
&& +\frac{(1-\beta)t}{(1-\alpha\eta)-\beta(1-\eta)} \int_0^1(1-s)y(s)ds
-\frac{(\alpha-\beta)t} {(1-\alpha\eta)-\beta(1-\eta)} \int_0^\eta (\eta-s)y(s)ds\nonumber\\
&&-\int_0^t (t-s)y(s)ds\nonumber\\
&=&-\int_0^t (t-s)y(s)ds+\frac{(\beta-\alpha)t-\beta }{(1-\alpha\eta)-\beta(1-\eta)} \int_0^\eta (\eta-s)y(s)ds\nonumber\\
&&+\frac{(1-\beta)t+\beta\eta}{(1-\alpha\eta)-\beta(1-\eta)}\int_0^1 (1-s)y(s)ds\nonumber\\
&=&\int_0^1 g(t,s)y(s)ds+\frac{\beta+(\alpha-\beta)t }{(1-\alpha\eta)-\beta(1-\eta)} \int_0^1 g(\eta,s)y(s)ds.\nonumber 
\end{eqnarray}
Then, $u(t)=\int_0^1 G(t,s)y(s)ds.$
The function $u$ presented above is the unique solution to the
problem (\ref{e2.1}), (\ref{e2.2}).}
\begin{lemma}\label{lemm2}\ \\
Let $0<\alpha<\frac{1}{\eta}$ and $0\le
\beta<\frac{1-\alpha\eta}{1-\eta}$.  If $y\in C([0, 1],\
[0, \infty))$, then the unique solution $u$ of the problem
(\ref{e2.1}), (\ref{e2.2}) satisfies
$$ u(t)\ge 0,\quad t\in[0, 1].$$
\end{lemma}
\proof{We know that if $u''(t)=-y(t)\le 0$ for $t\in (0,1)$, $ u(0)\ge 0$ and $ u(1)\ge 0,$  then $ u(t)\ge 0$ for $t\in [0,1].$ We have
$$\begin{array}{ccl}
u(0)&=&\displaystyle\frac{-\beta }{(1-\alpha\eta)-\beta(1-\eta)} \int_0^\eta (\eta-s)y(s)ds
+\frac{\beta\eta}{(1-\alpha\eta)-\beta(1-\eta)}\int_0^1 (1-s)y(s)ds
\\
&=&\displaystyle\frac{\beta }{(1-\alpha\eta)-\beta(1-\eta)}[ -\int_0^\eta (\eta-s)y(s)ds
+\eta\int_0^\eta (1-s)y(s)ds]\\
&&\displaystyle+\frac{\beta\eta}{(1-\alpha\eta)-\beta(1-\eta)}\int_\eta^1 (1-s)y(s)ds
\\
&=&\displaystyle\frac{\beta}{(1-\alpha\eta)-\beta(1-\eta)}[ \int_0^\eta s(1-\eta)y(s)ds]+\frac{\beta\eta}{(1-\alpha\eta)-\beta(1-\eta)}\int_\eta^1 (1-s)y(s)ds\geq0
\end{array}$$
and
$$\begin{array}{ccl}
u(1)&=&\displaystyle-\int_0^1 (1-s)y(s)ds+\frac{(\beta-\alpha)-\beta }{(1-\alpha\eta)-\beta(1-\eta)} \int_0^\eta (\eta-s)y(s)ds\\
&&\displaystyle+\frac{(1-\beta)+\beta\eta}{(1-\alpha\eta)-\beta(1-\eta)}\int_0^1 (1-s)y(s)ds
\\
&=&\displaystyle\frac{\alpha}{(1-\alpha\eta)-\beta(1-\eta)}[\eta\int_0^1 (1-s)y(s)ds+\int_0^\eta (\eta-s)y(s)ds]
\\
&\geq&\displaystyle\frac{\alpha}{(1-\alpha\eta)-\beta(1-\eta)}[\eta\int_0^\eta (1-s)y(s)ds+\int_0^\eta (\eta-s)y(s)ds]
\\
&=&\displaystyle\frac{\alpha}{(1-\alpha\eta)-\beta(1-\eta)}\int_0^\eta s(1-\eta)y(s)ds\geq 0.
\end{array}$$
Then, $u(t)\ge 0$ $\forall t\in[0, 1].$}
\begin{lemma} \label{l3}
The function $g$ has the following properties
\begin{description}
\item[$(i)$] $0 \leq g(t,s) \leq s(1-s)=g(s,s) \quad \forall t,s \in
[0,1].$ 
\item[$(ii)$]  Let $\theta\in [0,\frac{1}{2}]$. Then, for $t\in [\theta, 1-\theta]$ and $s\in[0,1],$ we have $$ g(t,s) \geq \min\{t,1-t\}g(s,s)\geq\theta g(s,s).$$
\end{description}
\end{lemma}
\proof{For $0\leq s \leq t \leq 1,$ we have
$$0\leq g(t,s)  =\displaystyle s(1-t) \leq s(1-s)=g(s,s).$$  And for $0\leq t \leq  s \leq 1,$ we have
$$g(t,s)  = t (1-s)\leq s (1-s)=g(s,s).$$
\noindent Thus $(i)$ holds.\\
\noindent If $s=0$ or $s=1$, we show that $(ii)$ holds.\\ For $0< s \leq t  \leq 1$ and $s\ne 1$ we have
$$\frac{g(t,s)}{g(s,s)} =    \displaystyle   \frac{t(1-s)}{s(1-s)}=\displaystyle\frac{t}{s}\geq  \displaystyle t \quad \forall
t \in [0,1].$$ \\ For  $\displaystyle 0\leq t \leq s < 1$ and $s\ne 0$ we have
$$\displaystyle\frac{g(t,s)}{g(s,s)}  =\displaystyle \frac{s(1-t)}{s(1-s)}=\displaystyle \frac{(1-t)}{(1-s)}  \geq   \displaystyle(1-t)\quad
  \forall t \in [0,1].$$ Then $$ g(t,s) \geq \min\{t,1-t\}g(s,s).$$ Thus, there exist $\theta \in]0,\frac{1}{2}]$ such that $$\displaystyle\frac{g(t,s)}{g(s,s)}\geq \theta,\quad  \forall t\in [\theta,1-\theta]$$
\noindent Thus $(ii)$ holds.}
\begin{lemma} \label{l4}
The function $G$ has the following properties
\begin{description}
\item[$(i)$] $G(t,s) \geq 0 \quad \forall t,s \in [0,1], $
\item[$(ii)$] $G(t,s) \leq k_1 g(s,s) \quad \forall t,s \in [0,1] $
and $k_1 = \displaystyle 1+\frac{\max\{\alpha,\beta\} }{\displaystyle(1-\alpha \eta)-\beta(1-\eta) },$
\item[$(iii)$] $\displaystyle \min_{\theta\leq t \leq 1-\theta} G(t,s) \geq
k_2 g(s,s) \quad \forall t,s \in [0,1] $ where $\theta\in (0,\frac{1}{2})$ and\\
$\displaystyle k_2 = \theta\left[1+\frac{\beta+\min\{(\alpha-\beta)\theta,(\alpha-\beta)(1-\theta)\}} {(1-\alpha\eta)-\beta(1-\eta)}\right]$
\end{description}
\end{lemma}
\proof{\\
$(i)$ From equation (\ref{e2.3}) and $(i)$ of Lemma \ref{l3}, we get
$$G(t,s)\geq 0\ \forall t,s\in [0,1].$$
$(ii)$ By equation (\ref{e2.3}) and $(i)$ of Lemma \ref{l3}, we have
$$G(t,s)= g(t,s)+\displaystyle \frac{\beta+(\alpha-\beta)t}{\displaystyle(1-\alpha \eta)-\beta(1-\eta)} g(\eta,s)$$
$$\leq g(s,s)+\displaystyle \frac{\max(\alpha,\beta)}{\displaystyle(1-\alpha \eta)-\beta(1-\eta)} g(s,s) = k_1 g(s,s).$$
$(iii)$ From $(ii)$ of Lemma \ref{l3}, for $t\in [\theta, 1-\theta]$  we have
$$\begin{array}{ccl}
\displaystyle G(t,s)
 & = &  \displaystyle g(t,s)+\displaystyle \frac{\beta+(\alpha-\beta)t}{\displaystyle(1-\alpha \eta)-\beta(1-\eta)} g(\eta,s)\\
 &\geq&  \theta g(s,s)+\displaystyle \frac{\beta+\min\{(\alpha-\beta)\theta,(\alpha-\beta)(1-\theta)\}}{\displaystyle(1-\alpha \eta)-\beta(1-\eta)}\theta g(s,s) \\
 &\geq&  \displaystyle  \theta\left [1+\displaystyle \frac{\beta+\min\{(\alpha-\beta)\theta,(\alpha-\beta)(1-\theta)\}}{\displaystyle(1-\alpha \eta)-\beta(1-\eta)} \right ]g(s,s)  =   k_2 g(s,s).
\end{array}
$$
}
\begin{lemma}\label{l5}
If $y \in C([0,1])$ and $y\geq 0$, then the unique solution $u$ of
the boundary value problem (\ref{e2.1}), (\ref{e2.2}) satisfies
$\displaystyle \min_{\theta\leq t \leq 1-\theta} u(t) \geq \gamma \|u\|_{1}$
where $\|u\|_{1}:= \displaystyle\sup\{|u(t)|;\ 0\leq t\leq 1\}$ and $\gamma := \displaystyle \frac{k_2}{k_1}.$
\end{lemma}
\proof{For any $t \in [0,1],$ by Lemma
\ref{l4} we have
\[u(t)  =  \displaystyle \int_0^1 G(t,s) y(s) ds \ \\
 \leq  \displaystyle k_1 \int_0^1 g(s,s) y(s) ds,
\  \]
\noindent thus $||u||_{1}  \leq k_1 \int_0^1 g(s,s) y(s) ds.$
\noindent Moreover, from $(iii)$ of Lemma \ref{l4} for $t \in [\theta,1-\theta],$  we have
\[u(t)  =  \displaystyle \int_0^1 G(t,s) y(s) ds \
\geq  \displaystyle k_2 \int_0^{1} g(s,s) y(s) ds \
\geq \displaystyle \frac{k_2}{k_1} ||u||_{1}.
\] Therefore $\displaystyle \min_{\theta\leq t \leq 1-\theta} u(t) \geq \gamma \|u\|_{1}.$\\ By Lemma \ref{lemm1}, we can show that the BVP (\ref{e2.1}), (\ref{e2.2}) has a solution
$u = u(t)$ if and only if u is a solution of the operator equation $u = Tu,$
where
$$
Tu(t)= \left\{
\begin{array}{ll}
\beta u(\eta), & -\tau \leq t \leq 0 , \\
\displaystyle \lambda\int_0^1 G(t,s) a(s)  f(s, u(s-\tau)) ds,& 0\leq  t \leq 1.
\end{array}
\right.
$$}
We assume the following hypothesis:
\begin{enumerate}
\item[$(H_{1})$] $f\in C([0,1]\times[0,\infty); [0,\infty))$,
\item[$(H_{2})$] $a\in C([0,1];[0,\infty))$ and there exists $t_0\in (0,1)$
such that $a(t_0)>0,$
\end{enumerate}
Let define,
$$f^{0}:=\displaystyle \limsup_{u \to 0}\max_{t \in [0,1]}\frac{f(t,u)}{u},\ \ \ \ f^{\infty}:=\displaystyle \limsup_{u \to \infty}\max_{t \in [0,1]}\frac{f(t,u)}{u},$$
$$ M_{1}: =\beta\displaystyle\int_{0}^{\tau}g(s,s)a(s)ds+\int_{\tau}^{1}g(s,s) a(s) ds \mbox{ and } \  M_{2}: = \displaystyle  \int_{0}^{1}g(s,s) a(s)ds.$$
The proof of our main results is based upon an application of the following Leray-Schauder fixed point theorem.
\begin{theorem}\ ([\ref{y4}])\label{wx}\\
Let $\Omega$ be a convex subset of a Banach space $X$, $0\in \Omega$ and $\Phi : \Omega \rightarrow \Omega$ be a completely continuous operator.
Then either
\begin{enumerate}
\item $\Phi$ has at least one fixed point in $\Omega,$ or
\item the set $\{x\in \Omega / x=\mu \Phi x,\ 0<\mu<1 \}$ is unbounded.
\end{enumerate}
\end{theorem}
\section{Main results}
Let $X=C[-\tau, 1]$ be a Banach space with norm $||u|| = sup\{|u(t)| : -\tau \leq t \leq 1\}$.
\begin{theorem}\ \label{yx}\\
Assume $(H_{1})$ and $(H_{2})$ hold. If $f^{0}< \infty$, then the boundary value problem (\ref{wq}) has at least one positive solution.
\end{theorem}
\proof{
Choose $\epsilon >0$ such that
$(f^{0}+\epsilon) \lambda k_{1} M_1\leq 1$. Since $f^{0}< \infty$, then there exists constant $B > 0,$ such that $f(s,u) < (f^{0}+\epsilon) u$ for $0<u\leq B$.\\
Let $$ \Omega = \displaystyle\{ u\ / u \in C([-\tau,1]), u\geq 0 , \|u\| \leq B, \min_{\theta \leq t \leq 1-\theta}u(t) \geq \gamma \|u\|  \}.$$
Then $\Omega$ is a convex subset of $X.$ \\
For $u\in \Omega$, by Lemmas \ref{lemm2} and \ref{l5}, we know that
$Tu(t)\geq 0$ and $\displaystyle\min_{\theta \leq t \leq 1-\theta}(Tu)(t) \geq \gamma \|Tu\|$.\\ Moreover,
$$\begin{array}{ccl}
Tu &\leq&\displaystyle\lambda k_{1}\int_0^1 g(s,s)a(s)f(s,u(s-\tau))ds
\\
&\leq&\displaystyle\lambda(f^0+\epsilon) k_{1}\int_0^1 g(s,s)a(s)u(s-\tau)
\\
&=&\displaystyle\lambda(f^0+\epsilon) k_{1} \left(\int_0^\tau g(s,s)a(s)\beta u(\eta)ds+\int_\tau^1 g(s,s)a(s)u(s-\tau)ds\right)
\end{array}$$
$$\begin{array}{ccl}
&\leq& \displaystyle\lambda(f^0+\epsilon) k_{1}\left(\beta \int_0^\tau g(s,s)a(s)ds+\int_\tau^1 g(s,s) a(s)ds\right)\|u\|
\\
&\leq& \|u\|\leq B.
\end{array}$$
Thus, $\|Tu\|\leq B.$ Hence, $T\Omega\subset\Omega.$\\
We shall show that $T$ is completely continuous.\\ Suppose $u_{n}\rightarrow u\ (n\rightarrow \infty)$ and $u_n \in \mathrm{\Omega} \ \forall n\in \mathbf{N},$ then there exists $M>0$ such that $\|u_n\|\leq M.$\\ Since $f$ is continuous on $[0,1]\times[0,M],$ it is uniformly continuous.\\ Therefore, $\forall \varepsilon >0$ there exists  $\delta>0$ such that $|x-y|<\delta$ implies $|f(s,x)-f(s,y)|<\epsilon \ \forall s\in [0,1],\ x, y\in [0,M]$ and there exists  $N$ such that $\|u_n -u\|<\delta$ for $n>N,$ so $|f(s,u_{n}(s-\tau))-f(s,u(s-\tau))|<\varepsilon,$ for $n>N$ and $s\in[0,1].$\\ This implies
$$\begin{array}{ccl}
|Tu_n(t)-Tu(t)|&\leq& \lambda k_{1}\displaystyle \int_0^1 g(s,s) a(s) |f(s,u_n(s-\tau))- f(s, u(s-\tau)| ds
\\
&\leq&  \displaystyle \lambda\epsilon k_{1} \int_0^1 g(s,s) a(s)ds.
\end{array}$$
Therefore $T$ is continuous.\\
\noindent Let $D$ be any bounded subset of $\mathrm{\Omega}$, then
there exists $\gamma>0$ such that $||u|| \leq \gamma$ for all $u \in D.$\\
Since $f$ is continuous on $[0,1]\times[0,\gamma]$ there exists $L>0$ such that $|f(t,v)|<L \ \forall (t,v)\in[0,1]\times [0,\gamma].$\\ Consequently, for all $u\in D$ and $t\in [0.1]$ we have
$$
\begin{array}{ccl}
|Tu(t)|& \leq & \left | \displaystyle \lambda k_{1}\int_0^1 g(s,s) a(s) f(s, u(s-\tau)) ds  \right|\\
& \leq & \lambda k_{1} L \int_0^1 g(s,s)  a(s)  ds.
\end{array}
$$
Which implies the boundedness of $TD.$\\
Since $G$ is continuous on $[0,1]\times[0,1],$ it is uniformly continuous.\\ Then  $\forall \epsilon>0$ there exists $\delta>0$ such that $|t_1-t_2|<\delta$ implies that $|G(t_1,s)-G(t_2,s)|<\epsilon \ \forall s\in [0,1].$ So, if $u\in D,$
$\displaystyle |Tu(t_1)-Tu(t_2)|\leq \lambda\int_0^1 |G(t_1,s)-G(t_2,s)| a(s) f(s,u_{n}(s-\tau)) ds\leq  \lambda L\epsilon\int_0^1 a(s) ds.$\\
From the arbitrariness of  $\epsilon,$ we get the equicontinuity of $TD$.\\
The operator $T$ is completely continuous by the mean of the Ascoli-Arzela theorem.\\
For $u\in \Omega$ and $u=\mu Tu,$ $0<\mu<1,$ we have $u(t)=\mu Tu(t)< Tu(t)<B,$ which implies $\|u\|\leq B.$ So, $\{x\in \Omega/ x=\mu\Phi x,\ 0<\mu<1 \}$ is bounded.\\ By theorem \ref{wx}, we deduce that operator $T$ has at least one fixed point in $\Omega$.
 Thus the boundary value problem (\ref{wq}) has at least one positive solution.
}
\begin{remark}\ 
The conditions of Theorem \ref{yx} are weaker than those of Theorem 3.1 in [\ref{y3}].
\end{remark}
\begin{theorem}\ \label{xy} \\
Assume $(H_{1})-(H_{2})$ hold. If $f^{\infty}<\infty$ is satisfied, then the boundary value problem (\ref{wq}) has at least one positive solution.
\end{theorem}
\proof{
Choose $\epsilon >0$ such that
$(f^{\infty}+\epsilon) \lambda k_{1} M_1\leq \displaystyle\frac{1}{2}$. Since $f^{\infty}< \infty$, then there exists constant $N>0$, such that $f(s,u)<(f^{\infty}+\epsilon) u$ for $u>N$.\\
 Let $B>0$ such that
$$B\geq N+1+\displaystyle 2 \lambda k_{1} M_{2}\displaystyle\max_{\begin{array}{c}
          0\leq s\leq 1 \\ 0\leq u\leq N
        \end{array}}f(s,u).$$
Let $$\Omega = \{u/u\in C[-\tau ,1],\ u\geq 0,\ \|u\| \leq B,\displaystyle  \min_{\theta\leq t\leq 1-\theta} u(t) \geq\gamma \|u\| \}.$$
Then $\Omega$ is a convex subset of $X.$ \\
For $u\in \Omega$, by Lemmas \ref{lemm2} and \ref{l5}, we have
$Tu(t)\geq 0$ and $\displaystyle\min_{\theta \leq t \leq 1-\theta}(Tu)(t) \geq \gamma \|Tu\|$.\\ Moreover, for $u\in\Omega$, we have
$$\begin{array}{ccl}
Tu(t)&=&\displaystyle \lambda\int_0^1 G(t,s) a(s) f(s, u(s-\tau)) ds\\
&\leq& \displaystyle   \lambda k_{1}\int_0^1 g(s,s) a(s) \ f(s, u(s-\tau)) ds\\
&=& \displaystyle \lambda k_{1}\left(\int_{J_1 = \{ s \in [0,1] / u > N  \} }g(s,s) a(s) f(s, u(s-\tau)) ds\right.\\
&& + \displaystyle\left.\int_{J_2= \{s \in [0,1] / u \leq  N \}} g(s,s) a(s) f(s, u(s-\tau)) ds\right)\\
&\leq & \displaystyle  \lambda k_{1} \left(  \int_0^1 g(s,s) a(s) (f^\infty+\epsilon) u(s-\tau) ds + \displaystyle  \int_0^1 g(s,s) a(s) \displaystyle\max_{\begin{array}{c}0\leq s\leq 1 \\ 0\leq u\leq N\end{array}}f(s, u(s-\tau)) ds\right)
\\
&\leq & \displaystyle\lambda k_{1}\left((f^\infty+\epsilon) \left[ \beta  \int_0^{\tau} g(s,s) a(s)   ds + \int_{\tau}^1 g(s,s) a(s)   ds \right] \|u \|\right.
\\
&&\left.\displaystyle+   \int_0^1 g(s,s) a(s) \displaystyle\max_{\begin{array}{c}0\leq s\leq 1 \\ 0\leq u\leq N\end{array}}f(s, u(s-\tau)) ds\right)
\end{array}$$
$$\begin{array}{ccl}
&\leq&   \lambda(f^\infty+\epsilon)  k_{1}  M_1 B +\lambda k_{1}  M_2  \displaystyle\max_{\begin{array}{c}0\leq s\leq 1 \\ 0\leq u\leq N\end{array}}  f(s, u(s-\tau)) \leq\frac{B}{2}  +  \frac{B}{2} = B.
\end{array}$$
Thus, $\|Tu\|\leq B.$ Hence, $T\Omega\subset\Omega.$\\
We can show that $T:\Omega\rightarrow\Omega$ is completely continuous.\\
For $u\in \Omega$ and $u=\mu Tu,$ $0<\mu<1,$ we have $u(t)= \mu Tu(t)< Tu(t)<B,$ which implies $\|u\|\leq B.$ So, $\{x\in \Omega/ x=\mu\Phi x,\ 0<\mu<1 \}$ is bounded.\\ By theorem \ref{wx}, we show that the operator $T$ has at least one fixed point in $\Omega$.\\
 Thus, the boundary value problem (\ref{wq}) has at least one positive solution.
}
\begin{remark}\
The conditions of Theorem \ref{xy} are weaker than those of Theorem 3.2 in [\ref{y3}].
\end{remark}

\end{document}